\documentclass[12pt]{article}
\usepackage{amscd, amssymb, latexsym, amsmath}

\newtheorem{thm}{Theorem}

\newtheorem{lem}{Lemma}

\newtheorem{con}{Conjecture}

\numberwithin{thm}{section} \numberwithin{cor}{section}
\numberwithin{pro}{section} \numberwithin{lem}{section}
\numberwithin{dfn}{section}
\numberwithin{rem}{section} \numberwithin{equation}{section}

\newcommand{\R}{\mathbb R}
\newcommand{\C}{\mathbb C}
\newcommand{\Z}{\mathbb Z}

\begin{document}
\title{Some recent developments in Lagrangian \\mean curvature flows}
\author{Mu-Tao Wang \footnote{The author is partially supported by National Science
Foundation Grant DMS 0605115 and an Alfred P. Sloan Research Fellowship.}}

\date{April 2, 2008}
\maketitle \centerline{email:\, mtwang@math.columbia.edu }
\begin{abstract}
We review some recent results on the mean curvature flows of Lagrangian submanifolds from the perspective of geometric partial differential equations. These include global existence and convergence results, characterizations of first-time singularities, and constructions of self-similar solutions.
\end{abstract}

\section{Introduction}
\subsection{The mean curvature flow}
A distinguished normal vector field called the \textit{mean
curvature vector field} exists on any submanifold
of a Riemannian manifold. It is characterized as the unique
direction along which the area or the volume the submanifold would
be decreased most effectively. A submanifold is \textit{minimal} if the mean curvature vector
vanishes at each point. The mean curvature flow is an
evolution process that moves a submanifold by its mean curvature
vector field. This turns out to be a nonlinear
parabolic system of partial differential equations for the
position functions of the submanifold. Brakke \cite{br} pioneered
the formulation of weak solutions of the mean curvature flow in the setting of geometric measure theory
in the seventies. In the eighties, Huisken and his coauthors \cite{hu1}  took up the geometric PDE approach
to study the mean curvature flow in Riemannian manifolds. There were also level set formulations of viscosity solutions for hypersurfaces by  Chen-Giga-Goto \cite{cgg} and independently by Evans-Spruck \cite{es1} later. Soon the mean
curvature flow joined the rank of most studied geometric evolution
equations like the Ricci flow and the harmonic map heat flows. Indeed, in many
aspects, its development paralleled that of the Ricci flow.

Suppose $M$ is a $m$ dimensional Riemannian manifold with a Riemannian metric $\langle \cdot, \cdot\rangle$. Let $\Sigma$ be another $n$ dimensional smooth manifold with $n<m$. An immersion of $\Sigma$ in $M$  is given by a mapping $F:\Sigma \rightarrow M$ so that the differential $dF$ has full rank at each point of $\Sigma$. In terms of local coordinates $x^1\cdots x^n$ on $\Sigma$, this is equivalent to the matrix $g_{ij}=\langle \frac{\partial F}{\partial x^i}, \frac{\partial F}{\partial x^j}\rangle$ being positive definite. Indeed, $g_{ij}$ defines a Riemannian metric which makes the image $F(\Sigma)$ a Riemannian submanifold. The second fundamental form of $F(\Sigma)$ is given by the tensor $A=(\nabla_{\frac{\partial F}{\partial x^i}}\frac{\partial F}{\partial x^j})^\perp$ where $\nabla$ denotes the Levi-Civita connection of $M$ and $\perp$ denotes the normal part of a vector in the bundle $TM|_{F(\Sigma)}$. The mean curvature vector $H$ is then the trace of the second fundamental form, i.e.

\[H=g^{ij}( \nabla_{\frac{\partial F}{\partial x^i}}\frac{\partial F}{\partial x^j})^\perp\] where $g^{ij}$ is the inverse of $g_{ij}$.

We recall the volume of $F(\Sigma)$ is calculated by

\[Vol(F(\Sigma))=\int_\Sigma \sqrt{\det{g_{ij}}} dx^1\wedge \cdots \wedge dx^n.\]

Suppose $F:\Sigma\times[0, \epsilon)\rightarrow M$ is a family of immersion so the variation field $\frac{\partial F}{\partial s}$ is a normal vector field along $F(\Sigma, s)$. The variation of the volume is then given by

\[\frac{d}{ds}Vol(F(\Sigma, s))=-\int_{\Sigma}\langle \frac{\partial F}{\partial s}, H\rangle \sqrt{\det{g_{ij}}} dx^1\wedge \cdots \wedge dx^n.\]

The mean curvature flow deforms a submanifold in the direction of
the mean curvature vector field $H$. Namely, a family  of immersion  $F:\Sigma\times[0, T)\rightarrow M$ is said to form a mean curvature flow if

\[(\frac{\partial F}{\partial t}(x,t))^\perp={H(F(x,t))}.\]

 The flow is a nonlinear weakly parabolic
system for $F$ and is invariant under reparametrization of $\Sigma$. Indeed, by coupling with a diffeomorphism of $\Sigma$, the flow can be made into a normal deformation, i.e. $\frac{\partial F}{\partial t}(x,t)$ is always in the normal direction. One can establish the short-time existence for any
smooth compact initial data. For a normal deformation, a simple calculation shows

\[\frac{d}{dt}\sqrt{\det{g_{ij}}(F(x,t))}=-|H(F(x,t))|^2\sqrt{\det{g_{ij}}(F(x,t))}, \] from which it follows that being a immersion is preserved by the mean curvature flow. Integrating this equality gives

\begin{equation}\label{volume}\frac{d}{dt} Vol(F(\Sigma, t))=-\int_{\Sigma} |H(F(x,t))|^2 \sqrt{\det{g_{ij}}(F(x,t))} dx^1\wedge \cdots \wedge dx^n.\end{equation}

In the rest of the paper, when there is no confusion, we will not differentiate between an immersed submanifold $\Sigma$ and its image $F(\Sigma)$. $T\Sigma$ will denote the tangent bundle and $N\Sigma$ will denote the normal bundle of $\Sigma$. We shall denote the image  $F(\Sigma, t)$ of a mean curvature flow by $\Sigma_t$ and the volume form $\sqrt{\det{g_{ij}}(F(x,t))} dx^1\wedge \cdots \wedge dx^n$ by $dv_{\Sigma_t}$.

\subsection{Co-dimension one vs. higher co-dimension}

The Lagrangian mean curvature flow is a special case of the mean curvature flow
in general. They are mostly higher co-dimensional in the sense that $m>n+1$ where $m$ is the dimension of the ambient manifold and $n$ is the dimension of the evolving submanifold.  There are abundant results of hypersurface (i.e. $m=n+1$) mean curvatures flows, while relatively little is known in the higher co-dimensional case.
Indeed, many techniques and results for hypersurface flows do not generalize
to higher co-dimensions.
The contrast between the hypersurface and higher co-dimensional
can be seen from two points. Firstly, the hypersurface case corresponds to a
scalar equation and the maximum principle holds in the following sense. Two embedded hypersurfaces evolving by the mean curvature flow will avoid each other and an embedded hypersurface remains embedded along the mean curvature flow. These are no longer true for higher co-dimensional mean curvature flows.
Secondly, the second fundamental form of a hypersurface is a
symmetric two tensor and various convexity condition associated with this
tensor has play an important role. However, for a higher
co-dimensional submanifold, the second fundamental is a symmetric
two tensor valued in the normal bundle and there is no natural
convexity condition for such a tensor. Of course, both these points are
related to the complexity of the normal bundle. The normal
bundle of an embedded orientable hypersurface is always trivial,
while the normal bundle of a higher co-dimensional submanifold could be
highly non-trivial.

The good news is that many great results like Brakke's regularity theorem \cite{br}, Hamilton's maximum principle for tensors \cite{ha}, Huisken's monotonicity formula \cite{hu2}, and White's regularity theorem \cite{wh} are valid in any dimension and co-dimension.

\subsection{What is special about being Lagrangian?}

A Lagrangian submanifold sits in a symplectic manifold $M$. Recall a symplectic manifold is a smooth manifold equipped with a closed non-degenerate two-form $\omega (\cdot, \cdot)$. A Lagrangian submanifold $\Sigma$ is characterized by the vanishing of $\omega|_{\Sigma}$.  To get the volume functional into the context, we consider on $M$ a Riemannian metric $\langle\cdot, \cdot\rangle$ and a compatible almost complex structure $J$ in the sense that $\omega(\cdot, \cdot)=\langle J(\cdot), \cdot \rangle$.

The simplest examples are Lagrangian submanifolds in the complex Euclidean space $\C^n$. Suppose
$z^i=x^i+y^i, i=1\cdots n$ are complex coordinates on $\C^n$. Let $\omega=\sum_{i=1}^n dx^i\wedge
dy^i$ be the standard symplectic form on $\C^n$ and
$J(\frac{\partial}{\partial x^i})=\frac{\partial}{\partial y^i}$
be the standard almost complex structure.  We have $\omega(X,
Y)=\langle JX, Y\rangle$ where $\langle \cdot, \cdot \rangle$ is the standard metric on $\C^n$. Recall a Lagrangian subspace of $\C^n$
is a subspace on which $\omega$ restricts to zero. As the form $\omega$ is invariant under the unitary group $U(n)$, the set of all Lagrangian subspaces in $\C^n$, so called the Lagrangian Grassmannian, is isomorphic to the homogeneous space $U(n)/SO(n)$.  Therefore a Lagrangian submanifold is simply a submanifold whose tangent spaces are all Lagrangian subspaces in $\C^n$.

Two prominent classes of Lagrangian submanifolds are the
followings.

\noindent 1). Graphs of symplectomorphisms: In this case, the ambient space is
the product of two symplectic manifolds. Take a smooth map
$f:(\C^n, \omega_1) \rightarrow (\C^n, \omega_2) $ such that $f^*\omega_2=\omega_1$ and the graph of $f$ in $(\C^n\times \C^n, \omega_1-\omega_2)$ is such an example.

\noindent 2). Graphs of one-forms: In this case, the ambient space is the
cotangent bundle of a symplectic manifold. Take a smooth map
$f:\R^n\rightarrow \R^n$ such that $f=\nabla u$ for some scalar function
$u:\R^n\rightarrow \R$. Identify the first $\R^n$ with the real part of $\C^n$ with coordinates $x^i$ and the second $\R^n$ as the imaginary part with coordinates $y^i$. The graph of $f$ in $(\C^n=\R^n\oplus \sqrt{-1}\R^n, \omega)$ with $\omega=\Sigma_{i=1}^n dx^i\wedge dy^i$ is such an example.

Two remarkable properties of Lagrangian mean curvature flow makes
one speculate that it perhaps behave in a better way than other
mean curvature flows in higher codimension. Firstly, the normal
bundle of a Lagrangian submanifold is canonically isometric to the
tangent bundle by the almost complex structure $J:T\Sigma\rightarrow  N\Sigma$. This gives a simpler description of the second fundamental form as a fully symmetric three tensor in $\odot^3 T\Sigma$, the coefficients of which being

\[h_{ijk}=\langle \nabla_{\frac{\partial F}{\partial x^i}}\frac{\partial F}{\partial x^j},  J(\frac{\partial F}{\partial x^k})\rangle.\]

As $JH$ becomes a tangent vector, it is dual to a one-form $\sigma$ on $\Sigma$; they are related by
$\sigma(\cdot)=\omega(H, \cdot)$.

 Secondly, when the ambient space is K\"ahler-Einstein, being Lagrangian is a condition that is
preserved along the mean curvature flow \cite{sm1}. In many occasions, one
is tempted to compare the Lagrangian mean curvature flow to the
K\"ahler-Ricci flow. This will be elaborated in the next section when we consider a Calabi-Yau ambient manifold.

\subsection{Calabi-Yau case}

The Lagrangian mean curvature flow is of particular interest when the ambient space is Calabi-Yau.
 Let $(M, g, \omega, J)$ be a real $2n$ dimensional
K\"ahler-Einstein manifold; i.e. the Ricci form is a multiple of the
the K\"ahler form:

\[Ric=c\omega.\]

Let $\Sigma$ be a $n$ dimensional Lagrangian submanifold of $M$.
 By the Codazzi
equation, we have
\[d\sigma=Ric|_{\Sigma}=c\omega|_{\Sigma}=0.\] Thus $\sigma$ is a closed one-form and
defines a cohomology class
$[\sigma]$.

When $M$ is Calabi-Yau with a canonical parallel holomorphic
$(n,0)$ form $\Omega$. By suitable normalization, the restriction
of $\Omega$ on $\Sigma$ gives a multi-valued function $\theta$,
called the {\it phase} function. Indeed, $*_\Sigma \Omega=e^{i\theta}$
where $*_\Sigma$ denote the Hodge star operator on $\Sigma$. It turns out the mean curvature form is
$\sigma=d\theta$. $[\sigma]$ is called the {\it Maslov class} and
can be defined through the Gauss map of $\Sigma$.

Lagrangian submanifolds with vanishing mean curvature $H$
are {\it minimal Lagrangian} submanifolds. When $M$ is Calabi-Yau,
a connected minimal Lagrangian has constant phase and is called a {\it
special Lagrangian} which corresponds to a class of calibrated
submanifolds first studied by Harvey and Lawson \cite{hl}.  They
also play important roles in the SYZ \cite{syz} conjecture in
Mirror symmetry. Roughly speaking, people expect special
Lagrangians behave much like holomorphic curves. It is thus
desirable to have a general method of constructing special (or
minimal) Lagrangian submanifolds. Schoen and Wolfson \cite{sw1}
studied the existence problem by the variational method. Another approach more related to algebraic geometry and symplectic topology has been investigated by Joyce.

Based on the duality between graded Lagrangian submanifolds and
stable vector bundles, Thomas and Yau \cite{ty} made the following
conjecture:

\begin{con}

Let $M$ be Calabi-Yau and $\Sigma$ be a compact embedded
Lagrangian submanifold with zero Maslov class, then the mean
curvature flow of $\Sigma$ exists for all time and converges
smoothly to a special Lagrangian submanifold in the Hamiltonian
isotopy class of $\Sigma$.

\end{con}

Suppose $\Sigma_t$ is the mean curvature flow of $\Sigma$. The phase function
$\theta$ on $\Sigma_t$ is evolved by the heat equation

\begin{equation}\label{phase_heat}\frac{d\theta}{dt}=\Delta_{\Sigma_t} \theta,\end{equation} where $\Delta_{\Sigma_t}$ is the Laplace operator of the induced
metric on $\Sigma_t$. Thus being of zero Maslov class is preserved along the flow.
This is a rather bold conjecture as it is easy to see that the
mean curvature flow of any compact submanifold of the Euclidean space develops
finite time singularities. The most common singularity is the
so-called \textit{neck-pinching}. Without the assumption on the Maslov class, Schoen and Wolfson \cite{sw2} construct example that develop such singularities in finite time. We shall come back to this point in \S 2.3.

\subsection{Overview of the article}
In this paper, we review some recent results on the mean
curvature flows of Lagrangian submanifolds.
The review is by no means comprehensive or complete and is subject to the author's personal preference. Many other interesting works on the Lagrangian mean curvature flow are not discussed in the article, e.g. \cite{gssz}, \cite{lw}
,\cite{ne2}, and \cite{pa}.

A fundamental question in geometric flows is under what conditions on the initial data can we prove the global
existence and convergence in the smooth category. The mean curvature flow, as a quasi-linear system, forms singularity exactly
when the second fundamental form of the submanifold blows up, i.e. $\Sigma_t$ becomes singular as
$t$ approaches $T$ if and only if $\lim_{t\rightarrow T}
\sup_{\Sigma_t} |A|^2 \rightarrow \infty$. In \S 2.1, we discuss global existence and convergence results for special initial data that correspond to the two classes of Lagrangian submanifolds discussed in \S 1.3.

  For more general initial data set, we need to understand singularity formations through the blow-up analysis. Suppose the flow exists on $[0, T)$ and $\sup_{\Sigma_t}|A|^2\rightarrow \infty$ as $t\rightarrow T$, we perform parabolic blow-up of the solution near $T$. For simplicity, here we demonstrate the idea in the case when the ambient space is the Euclidean space $\R^N$.
  This process depends on three parameters $t_i$(when), $y_i$(where), and  $\lambda_i$ (how much) and eventually we let $\lambda_i\rightarrow \infty$.  Take the space time track of the mean curvature flow $\frak{M}=\cup_{t\in [0, T)}\Sigma_t$ in $\R^N\times \R $ and consider the map
\[\R^N\times[0, T)\rightarrow \R^N\times [-\lambda_i^2 t_i,
\lambda_i^2(T-t_i))\]
by sending $(y,t)\mapsto (\lambda_i(y-y_i), \lambda_i^2(t-t_i))$ and thus $(y_i, t_i)\mapsto (0,0)$.

The image of $\Sigma_t$ can be described as
\[{\Sigma}^i_{s}=\lambda_i(\Sigma_{t_i+\frac{s}{\lambda_i^2}}
-y_i)\]
 in which $s=\lambda_i^2(t-t_i)$. It turns out the space-time track $\frak{M}_i=\cup_{s}{\Sigma}^i_{s}$ forms another mean curvature flow by the invariance of the scaling that lives in $[-\lambda_i^2 t_i,
\lambda_i^2(T-t_i))$.

In order to obtain smooth limit, two types of parabolic blow-ups are often used depending on how fast the second fundamental form blow up.

{\it Type I blow-up}: Also called a central blow-up where the center $(t_i, y_i)=(T, y_0)$ is fixed and $\lambda_i\rightarrow \infty$. When $|A|^2(T-t)$ is bounded, the limit is smooth.
Nevertheless, the limit always exists weakly in the sense of geometric measure theory and is an ancient self-similar solution that lives in $(-\infty, 0]$. This is the parabolic analogue of a cone.

{\it Type II blow up}: The blow-up center $(y_i, t_i)$ is at a point where $|A|^2(T-t)$ almost achieves its
maximum. The scale is proportional to $|A|^2$ so we get a smooth limit that often lives on $(-\infty, \infty)$ with uniformly bounded second fundamental form, so called an eternal solution. In \S 2.2, we discuss the characterization of first time singularity for Lagrangian mean curvature flow under the Type I blow-up procedure.

The singularities are often classified into \textit{type I singularity} or \textit{type II singularity} according to whether $\sup_{\Sigma_t}|A|^2(T-t)$ is bounded or unbounded as $t\rightarrow T$ (not to be confused with the  type I and type II blow-up procedures). The simplest Lagrangian mean curvature flow reduces to the so-called curve
shortening flow on a two-dimensional orientable surface as any curve is Lagrangian. In this case, it can be proved that any type I singularity is a shrinking circle and any type II singularity is a Grim Reaper, both are self-similar solutions. It is thus concievable that the eventual understandings of the singularities will rely on the classifications of self-similar solutions. In \S 2.3, we discuss the constructions of self-similar solutions in Lagrangian mean curvature flows.

The author would like to thank B. Andrews, K. Ecker, R. Hamilton, M. Haskin, G. Huisken, T. Ilmanen, D. Joyce, Y.-I. Lee, N. C. Leung, A. Neves, K. Smoczyk, M.-P. Tsui, T. Y. H. Wan, and B. White for helpful discussions on this subject.

\section{Results}
\subsection{Global existence and convergence}
Given an immersed submanifold $\Sigma_0$ of a Riemannian manifold $M$, we ask when we can find a family of immersions that forms a mean curvature flow $\Sigma_t$ and when $\Sigma_t \rightarrow \Sigma_\infty$ in $C^\infty$ for a smooth immersed submanifold $\Sigma_\infty$. As was remarked in the overview of the article, this boils down to bounding the second fundamental form for all $t\in [0, \infty)$ and as $t\rightarrow \infty$. The one-dimensional curve-shortening flow is a well-studied area and there are many beautiful global existence and convergence results by e.g. Gage-Hamilton \cite{gh} and Grayson\cite{gr1}\cite{gr2}. We refer to the book by Chou-Zhu \cite{cz} for results in this direction.

The next simplest case will be two-dimensional Lagrangian surfaces in a four-dimensional symplectic manifold. We recall that the graph of a symplectomorphism is naturally a Lagrangian submanifold of the product space. In this case, there is the following theorem:

\begin{thm}
Let $(\Sigma^{(1)}, \omega_1)$ and $(\Sigma^{(2)},\omega_2)$ be two
diffeomorphic compact Riemann surface of the same constant
curvature $c$. Suppose $\Sigma$ is the graph of a
symplectomorphism $f:\Sigma^{(1)}\rightarrow \Sigma^{(2)}$ as a Lagrangian
submanifold of $M=(\Sigma^{(1)}\times \Sigma^{(2)}, \omega=\omega_1-\omega_2)$
and $\Sigma_t$ is the mean curvature flow with initial surface
$\Sigma_0=\Sigma$. Then $\Sigma_t$ remains the graph of a
symplectomorphism $f_t$ along the mean curvature flow. The flow
exists smoothly for all time and $\Sigma_t$ converges smoothly to
a minimal Lagrangian submanifold as $t\rightarrow \infty$.
\end{thm}
The assumption on the curvatures of $\Sigma^{(1)}$ and $\Sigma^{(2)}$ makes $M$ a K\"ahler-Einstein manifold with the product metric.
The long time existence for all cases and the smooth convergence for $c>0$  was  proved in
\cite{wa2}. The smooth convergence for $c\leq 0 $ was established in \cite{wa4}

Assuming an extra angle condition, Smoczyk \cite{sm2} proved the theorem when $c\leq 0$. Smoczyk's proof of the global existence and convergence result is different from that of \cite{wa2} and $\cite{wa4}$. Instead of applying blow-up analysis, he proved directly by the maximum principle that the second fundamental form is uniformly bounded independent of time. This gives the global existence and the convergence at the same time.

We remark the existence of such minimal Lagrangian submanifold was
proved using variational method by Schoen \cite{sc} (see also Lee
\cite{le}).

In the following, we briefly describe the proof of the theorem. The proof is divided into three parts.

\noindent 1) $\Sigma_t$ remains the graph of a symplectomorphism $f_t$ as long as the flow exists smoothly.

Since $M$ is K\"ahler-Einstein, $\Sigma_t$ remains a Lagrangian surface. This can indeed be shown by considering the evolution equation satisfied by the function
$*\omega(p)=\omega (e_1, e_2)$ where $\{e_1, e_2\}$ is any
oriented orthonormal basis for $T_p\Sigma$. Likewise, we can consider the heat equation satisfied by the function
$*\omega_1(p)=\omega_1(e_1, e_2)$.   In an orthonormal basis, we can represent the second fundamental form $A$ by $h_{ijk}=\langle \nabla_{e_i} e_j, J(e_k)\rangle$ and the mean curvature vector by $H_k=\langle H, J(e_k)\rangle=\sum_{i=1}^2 h_{iik}$. It was computed in \cite{wa2} that $\eta_t=2*_{\Sigma_t}\omega_1$ satisfies

\begin{equation}\label{eta_t}\frac{d}{dt}\eta_t=\Delta
\eta_t+\eta_t(2|A|^2-|H|^2)+c\eta_t(1-\eta_t^2)\end{equation} where $|A|^2=\sum_{i,j,k}h_{ijk}^2$ and $|H|^2=\sum_{k}H_k^2$ are the squared norms of the second fundamental form and the mean curvature vector, respectively. As $|H|^2\leq 2|A|^2$, $\eta_t>0$ is preserved along the flow by the maximum
principle.

Notice that $*\omega_1$ is in fact the Jacobian of
$\pi_1|_\Sigma$ where $\pi_1:\Sigma^{(1)}\times \Sigma^{(2)}\rightarrow
\Sigma^{(1)}$ is the projection map onto the first factor. By the inverse function theorem, $*\eta_t>0$ if and only if
$\Sigma_t$ can be locally written a graph over $\Sigma_1$. Therefore $\eta_t>0$ implies $\Sigma_t$ is the graph of a symplectomorphism.

 Indeed, not only $\eta_t>0$ but by comparing to solutions of the ordinary differential equation $\frac{d}{dt} f=cf(1-f^2)$, we arrive at
\[\eta_t\geq \frac{\alpha e^ct}{\sqrt{1+\alpha^2e^{2ct}}}\] where
$\alpha$ is a constant that satisfies
$\frac{\alpha}{\sqrt{1+\alpha^2}}=\min_{\Sigma_0}\eta.$ In particular, $\eta_t\rightarrow 1$ as $t\rightarrow
\infty$ uniformly when $c>0$.

\noindent 2) Global existence for all finite time

Utilizing the full symmetry of the second fundamental form $h_{ijk}$, one can show that

\[|H|^2\leq \frac{4}{3} |A|^2.\]
Therefore,
\begin{equation}\label{eta}\frac{d}{dt}\eta_t\geq \Delta
\eta_t+\frac{2}{3}\eta_t|A|^2+c\eta_t(1-\eta_t^2).\end{equation}

We can apply the type I blow-up procedure to this solution at any space-time point.
Equation \ref{eta} and the positive lower bound of $\eta_t$ at any finite time will imply the integral of $|A|^2$ vanishes on the type-I blow-up limit. It follows from White's regularity theorem \cite{wh} that any such point is a regular point.

\noindent 3). Convergence at $t=\infty$.

  The aim is to bound $|A|^2$ as $t\rightarrow \infty$ since Simon's \cite{si} convergence theorem for gradient flows is applicable in this case. Suppose $\sup_{\Sigma_t}|A|^2\rightarrow \infty$, we apply the type II blow-up procedure to the solution at $t=\infty$. Pick a sequence of $t_i$ and point $p_i\in\Sigma_{t_i}$ such
 that the space-time track $\frak{M}_i$, after shifting $(p_i, t_i)$ to $(0, 0)$ and scaling by the factor
 $|A|(p_i, t_i)$, has uniformly bounded
 second fundamental form and $|A|(0,0)=1$.

 By compactness, $\frak{M}_i\rightarrow \frak{M}_\infty$, which is an eternal solution of the mean curvature flow defined on $(-\infty, \infty)$ with
 uniformly bounded second fundamental form and $|A|(0, 0)=1$. When $c>0$, recall $\eta_t\rightarrow 1$ as $t\rightarrow
 \infty$, this implies $\eta\equiv 1$ on the limit $\frak{M}_\infty$ and each time
 slice must be a flat space,  contradicting with $|A|(0,0)=1$.

 In case when $c\leq 0$, $\eta_t$ no longer converges to $1$ as $t\rightarrow \infty$ we consider instead the evolution equation for $|H|^2$. It was computed in \cite{wa4} that
\[(\frac{d}{dt}-\Delta)|H|^2=-2|\nabla
H|^2+2\sum_{ij}(\sum_k H_k h_{kij})^2+c(2-\eta^2)|H|^2.\]

Coupling with the equation for $\eta$ (\ref{eta_t}) and integrating over $\Sigma_t$ gives
 \[\frac{d}{dt}\int_{\Sigma_t} \frac{|H|^2}{\eta} dv_{\Sigma_t}\leq
c\int_{\Sigma_t}\frac{|H|^2}{\eta}dv_{\Sigma_t}.\]

This implies $\int_{\Sigma_t}|H|^2 dv_{\Sigma_t}\rightarrow 0$ as
$t\rightarrow \infty$ since  $\int_0^\infty \int_{\Sigma_t} |H| dv_{\Sigma_t} dt<\infty$. The limit $\frak{M}_\infty$ obtained
earlier has $\int |H|^2=0$ and thus each time slice is a minimal area
preserving map  from $\C$ to $\C$ which must be flat by a result of
Ni \cite{ni}.

In general dimension, Smoczyk and the author \cite{sw} proved a general global existence and convergence theorem for Lagrangian graphs in $T^{2n}$, a flat torus of dimension $2n$.

\begin{thm} Let $\Sigma$ be a
Lagrangian submanifold in $T^{2n}$. Suppose $\Sigma$ is the graph
of $f:T^n\rightarrow T^n$ and the potential function $u$ of $f$ is
convex. Then the mean curvature flow of
 $\Sigma$ exists for all time, remains a Lagrangian graph, and converges smoothly to a flat Lagrangian
submanifold.
\end{thm}

The flow in terms of the potential $u$ is a fully nonlinear parabolic equation:
\begin{equation}\label{spl}
\frac{du}{dt}=\frac{1}{\sqrt{-1}}\ln\frac{\det(I+\sqrt{-1}
D^2u)}{\sqrt{\det(I+(D^2u)^2)}}\end{equation} where $D^2 u$ is the Hessian of $u$. Notice $\frac{\det(I+\sqrt{-1} D^2u)}{\sqrt{\det(I+(D^2u)^2)}}$ is a unit complex number, so the right hand side is always real.
This theorem generalizes prior global existence and convergence results in general dimensions in \cite{wa3} and \cite{sm3}.

A important step in the proof is to show the convexity condition $D^2_{ij} u>0$ is preserved which we describe in the rest of this section. This involve interpreting the convexity condition as the positivity of some symmetric two tensor on $\Sigma_t$ and compute the parabolic equation with respect to the induced (evolving metric) on $\Sigma_t$. It turns out if we denote $\pi_1$ ($\pi_2$) to be the projection onto the first (second) factor of $T^n\times T^n$.
 The condition $D^2_{ij} u$ is the same as

\[\langle J\pi_1(X), \pi_2(X)\rangle >0\]for any $X\in T\Sigma$.
$S(\cdot, \cdot)=\langle J\pi_1(\cdot), \pi_2(\cdot)\rangle $ defines a two-tensor on $T^{2n}$ and the Lagrangian condition implies the restriction of $S$ to any Lagrangian submanifold is a symmetric tensor.

 \begin{lem} This positivity of $S$ is preserved along the
mean curvature flow, i.e. $S|_{\Sigma_t}>0$ for $t>0$ as long as $S|_{\Sigma_0}>0$.
\end{lem}

A direct approach is the calculate the evolution of $S|_{\Sigma_t}$ and apply Hamilton's maximal principle for tensors, see equation (3.3) in \cite{sw}.

Another more systematic approach is to study how the tangent space of $\Sigma_t$ evolves as this contains the information of $D^2_{ij}u$. Since the tangent space of $T^{2n}$ can be identified with $\C^n$. We may consider the Gauss map of $\Sigma_t$ given by
$\gamma_t:\Sigma_t\rightarrow LG(n)=U(n)\slash SO(n)$, the
Lagrangian Grassmannian, by sending  a point $p\in \Sigma$ to the tangent space $ T_p\Sigma\subset
\C^n$.
The following theorem is prove in \cite{wa5}
 \begin{thm}$\gamma_t$ is a harmonic map heat flow.
\end{thm}

 For an ordinary heat equation $\frac{d}{dt}f=\Delta f$, the conditions $f>0$ and
 $f=0$ are preserved by the maximum principle. For a harmonic map heat flow into a Riemannian manifolds, the analogy is that the image of the map will remain in a convex or totally geodesic set.

 Since $LG(n)$ is totally geodesic subset of the Grassmannian, this provides an alternate
way to show why being Lagrangian is preserved along the mean curvature flow. Also the determinant map from $U(n)\slash SO(n)$ to $U(1)\equiv S^1$ is totally geodesic. Thus the composition $\gamma_t\circ \det$ is a harmonic map heat flow into $S^1$. It is easy to see that $\gamma_t\circ \det$ is exactly the phase function $\theta$ and this is another way to derive (\ref{phase_heat}).

 To show $\{L\in LG(n), S|_L>0\}$ is a convex subset, we study the geodesic equation on $LG(n)$ and the details can be found in \cite{sw}. We remark that as being a minimal Lagrangian in $\C^n$ is an invariant property under the symmetry group $U(n)$. The equation of $u$ indeed enjoys more symmetric than a general fully non-linear Hessian equation. This observation provides more equivalent conditions under which the global existence and convergence theorems can be proved (see the last section in \cite{sw}).

\subsection{Characterization of first-time singularities}
In \cite{wa1}, the author introduced the notion of  \textit{almost
calibrated} Lagrangian submanifolds in the study of characterizing the first time singularity.   Recall for a special Lagrangian, the Lagrangian
angle (after a shifting) satisfies $\cos\theta=1$. A Lagrangian submanifolds in a
Calabi-Yau manifold is said to be ``almost calibrated" if
$\cos\theta\geq \epsilon $ for some $\epsilon>0$. This has proved to be a very useful condition in the
study of Lagrangian mean curvature flow.

As $\cos\theta$ satisfies

\begin{equation}\label{cos_phase}
\frac{d}{dt}\cos\theta=\Delta \cos\theta+\cos\theta |H|^2,
\end{equation}
being almost calibrated is another condition that
is preserved along the Lagrangian mean curvature flow. The following theorem is proved in \cite{wa1}.

\begin{thm}
An almost calibrated Lagrangian submanifold does not develop any type I singularity along the mean curvature flow.
\end{thm}

This is established by coupling equation (\ref{cos_phase}) with Huisken's monotonicity formula. In particular, no ``neck-pinching" will be forming in the Thomas-Yau conjecture if this condition is assumed.
To demonstrate the idea, let us pretend the Lagrangian submanifolds are compact and lie in $\R^n$. As along as characterizing finite time singularity is concerned, this does not pose any serious restriction as the ambient curvature will be scaled away under a blow-up procedure. Very mild assumption needs to be imposed at infinity to assure the integration by parts work.
Suppose the flow exists on $[0, \infty)$ and consider the backward heat kernel at $(y_0, T)$.
\begin{equation}
\rho_{y_0, T}(y,t)=\frac{1}{(4\pi(T-t))^{n\over 2}} \exp
(\frac{-|y-y_0|^2}{4(T-t)})
\end{equation}

Huisken's monotonicity formula implies
\[
\begin{split}
&\frac{d}{dt}\int_{\Sigma_t}  \rho_{y_0,T} dv_{\Sigma_t}=-\int_{\Sigma_t} |{H}+{1\over 2(T -t)} F^\perp|^2\rho_{y_0, T}dv_{\Sigma_t}.
\end{split}
\]

Coupling with the equation for $\cos\theta$ (\ref{cos_phase}), we obtain
\[
\begin{split}
&\frac{d}{dt}\int_{\Sigma_t} \rho_{y_0,T} (1-\cos\theta) dv_{\Sigma_t}\,
\\
=&-\int_{\Sigma_t} \rho_{y_0,T} |H|^2 \cos\theta dv_{\Sigma_t}-\int_{\Sigma_t}  |{H}+{1\over 2(T -t)} F^\perp|^2(1-
\cos\theta)dv_{\Sigma_t}.
\end{split}
\]

Both these equations are scaling invariant and continue to hold for the type I blow-up at $(y_0, T)$ (notice that $\theta$ is scaling invariant) and thus
 \begin{equation}\label{scale1}
 \begin{split}
\frac{d}{ds}\int_{\Sigma_s^i}  \rho =-\int_{\Sigma_s^i} \rho |{H}-{1\over 2s} F^\perp|^2
\end{split}
\end{equation}
\begin{equation}\label{scale2}\begin{split}
&\frac{d}{ds}\int_{\Sigma_s^i} \rho (1-\cos\theta)\,
\\=&-\int_{\Sigma_s^i} \rho |H|^2 \cos\theta-\int_{\Sigma_s^i}  \rho |{H}-{1\over 2s} F^\perp|^2(1-
\cos\theta)
\end{split}
\end{equation} where $\rho$ is the backward heat kernel at $(0,0)$ and $s=-\lambda_i^2(T-t)$. Therefore it is not hard to see that there exists a sequence of rescaled submanifold on which both
$\int |H|^2$ and  $\int  |{H}-{1\over 2s} F^\perp|^2$ are approaching zero. Thus $H=0$ and $F^\perp=0$ weakly on each time slice of the limit. This indicates that each time slice of the limit should be a union of minimal Lagrangian cones. If we assume
the singularity is of type I, then the limit flow is smooth and thus must be a flat space. White's regularity theorem implies the point is a smooth point.

Notice that (\ref{scale1}) implies local area bound and (\ref{scale2}) implies
local bound for the $L^2$ norm of mean curvatures on $\Sigma^i_s$. It follows from compactness theorems in geometric measure theory that the limit is rectifiable and this was carried out by Chen and Li in \cite{cl}.  As $|H|^2=|\nabla\theta|^2$, a natural question arises whether the phase $\theta$ is a constant on this union of minimal Lagrangian cone. Notice that even a union of special Lagrangian cones may have different phases and hence not necessarily area-minimizing.
Chen and Li \cite{cl} claimed that the phase function is a constant on the limit by proving a Poinc\'are inequality for $\theta$. Unfortunately, the proof of Theorem 5.1 in \cite{cl} overlooks some technical difficulties. Neves later give a different proof assuming two extra conditions and using the evolution equation of the Liouville form  $\lambda=\sum_{i=1}^n x^i dy^i-y^idx^i$. We refer to his paper \cite{ne1} for the precise statement of the theorem (Theorem B). Neves \cite{ne1} was also able to replace the assumption of almost calibrated by zero Maslov class by observing that the equation for $\cos\theta$ can be replaced by
\[\frac{\partial}{\partial t} \theta^2=\Delta \theta^2-2|H|^2.\]

\subsection{Constructions of self-similar solutions}

A important tool in the study of geometric flows is the blow-up analysis. A blow-up solution of the mean curvature flow sits in the Euclidean space and often enjoys more symmetry. It is important to study these special solutions as singularity models. A mean curvature flow in the Euclidean space is said to be \textit{self-similar} if it is moved by an ambient symmetry. We may consider ansatz of the type

\begin{equation}\label{scale}F(x, t)=\phi(t) F(x)\end{equation}

and

\begin{equation} \label{translate}F(x,t)=F(x)+\psi (t)\end{equation} which correspond to scaling symmetry and translating symmetry, respectively. The ansatz coupled with the mean curvature flow equation gives an elliptic equation for $F(x)$.  For solution of the form (\ref{scale}), $F(x, t)$ is called an expanding  or a shrinking soliton depending on whether $\phi(t)$ is greater or less than one, respectively.  A mean curvature flow $F(x, t)$ that satisfies (\ref{translate}) is called a translating soliton.

Henri Anciaux constructed examples of  Lagrangian shrinking and expanding solitons in \cite{an}. All the examples are based on minimal Legendrian immersions in $S^{2n-1}$ and
the solutions are asymptotic to the associated minimal Lagrangian
cones.

Yng-Ing Lee and the author \cite{lwa} constructed examples of self-similar shrinking and expanding
Lagrangian mean curvature flows that are asymptotic to Hamiltonian
stationary cones. They were able to glue them together to form weak solutions of the mean curvature flow in the sense of Brakke. In a new preprint of Joyce, Lee and Tsui \cite{jlt} constructed new examples of self-similar solutions, in particular translating solitons. \cite{nt} gave some characterizations of translating solitons in the two dimensional case.

\section{Prospects}

There have been several attempts to find counterexamples of the Thomas-Yau conjecture.
 Other than the examples of Schoen and Wolfson in \cite{sw2}, Neves \cite{ne1} constructed almost calibrated complete non-compact Lagrangian surfaces in $\C^2$ that develop finite time singularities. However, there is still no genuine counterexample to the Thomas-Yau conjecture as it was stated.

It should be noted that Schoen and Wolfson \cite{sw3} proved the following
existence result of special Lagrangians in a K-3 surface.

\begin{thm}
Let $X$ be a K3 surface with a Calabi-Yau metric. Suppose that $\gamma \in H_2(X; \Z)$ is a
Lagrangian class that can be represented by an embedded Lagrangian
surface. Then $\gamma$ can be represented by a special Lagrangian
surface.
\end{thm}

A mean curvature flow proof
of this theorem will confirm the Thomas-Yau
conjecture in two dimension. Since it was already shown that there is no type-I singularity, we need to focus on type II singularities in the zero Maslov class or almost calibrated case.  A general type-II singularity can be scaled to get an eternal solution with uniformly bounded second fundamental form that exists on $(-\infty, \infty)$. Such a solution of parabolic equation should be rather special and we hope to say more about it in the near future.

For a general initial data, one is tempted to speculate that, just as in the Ricci flow case, surgeries are necessary in order to continue the flow. It was commented in Perelman's paper \cite{pe} that when the surgery scale goes to zero, the solution with surgeries should converge to a ``weak
solution" of the Ricci flow, a notion that has yet to established. Weak formulations for the mean curvature
flow are available. However, as weak solutions are no longer unique, it is necessary to instruct the flow how to continue after singuarities. The examples found in Lee-Wang \cite{lwa} and Joyce-Lee-Tsui \cite{jlt} start out as
shrinking solitons as $t<0$, approach to Schoen-Wolfson cones as $t\rightarrow 0$ and resolve to expanding solitons for $t>0$. They altogether form a
Brakke flow. Notice that the Schoen-Wolfson cones are the only obstructions to the existence of regular minimal Lagrangian in two dimensional. It would be of great interest to glue in these models whenever such singularities form. We believed such models will play important roles in the global existence of surgery or weak solutions of the Lagrangian mean curvature flow.

\end{document}